\long\def\ignore#1\endignore{XXX}
\font\tenmsa=msam10
\font\tenmsb=msbm10

\font\largebf=cmbx10 scaled\magstep2
\font\largemi=cmmi10 scaled\magstep2

\def\all{\hbox{for all}}
\def\All{\hbox{For all}}
\def\and{\hbox{and}}

\def\bra#1#2{\langle#1,#2\rangle}
\def\Bra#1#2{\big\langle#1,#2\big\rangle}
\def\cite#1\endcite{[#1]}

\def\dbs{^{**}}

\def\eps{\varepsilon}
\def\exs{\hbox{there exists}}
\def\f#1#2{{#1 \over #2}}

\def\half{{\textstyle\f12}}

\def\infn{\inf\nolimits}
\def\Lt{{\wt L}}
\def\Ltt{L^{\rm TT}}

\def\lr{\Longrightarrow}

\def\on{\hbox{on}}

\def\phi{\varphi}

\def\pt{\wt p}
\def\qed{\hfill\hbox{\tenmsa\char03}}

\def\qLt{q_{\wt L}}
\def\quand{\quad\and\quad}
\def\r{\hbox{\tenmsb R}}
\def\rbar{\,]{-}\infty,\infty]}

\long\def\slant#1\endslant{{\sl#1}}
\def\st{\hbox{such that}}

\def\supn{\sup\nolimits}

\def\wh{\widehat}
\def\wt{\widetilde}

\def\defSection#1{}
\def\defCorollary#1{}
\def\defDefinition#1{}
\def\defExample#1{}
\def\defLemma#1{}
\def\defNotation#1{}
\def\defProblem#1{}
\def\defRemark#1{}
\def\defTheorem#1{}
\def\locno#1{}
\def\meqno#1{\eqno(#1)}
\def\nmbr#1{}
\def\Proof{\medbreak\noindent{\bf Proof.}\enspace}
\def\Proo{{\bf Proof.}\enspace}
\def\Signoff{}
\def \INTsec{1}
\def \SNLsec{2}
\def \PCdef{2.1}
\def \NORMdef{2.2}
\def \IOTAone{1}
\def \Pone{2}
\def \QCONTone{3}
\def \QMAXone{4}
\def \DUALone{5}
\def \ITone{6}
\def \QTone{7}
\def \DUALthree{8}
\def \Hex{2.3}
\def \MONNIsec{3}
\def \EEex{3.1}
\def \PTtwo{9}
\def \PTrem{3.2}
\def \NIdef{3.3}
\def \RECASTlem{3.4}
\def \CONVsec{4}
\def \Jone{10}
\def \Jtwo{11}
\def \LINsec{5}
\def \COERCIVITY{12}
\def \QCSTARlem{5.1}
\def \QClem{5.2}
\def \INFthm{5.3}
\def \INFone{13}
\def \INFtwo{14}
\def \INFthree{15}
\def \INFfour{16}
\def \INFcor{5.4}
\def \INFCzero{17}
\def \INFCone{18}
\def \INFCtwo{19}
\def \INFCthree{20}
\def \INFCfour{21}
\def \NEGcor{5.5}
\def \LINMONsec{6}
\def \RHOone{22}
\def \ADJPOLone{23}
\def \MONcor{6.1}
\def \MONone{24}
\def \MONfour{25}
\def \SURBBsec{7}
\def \SURone{26}
\def \SURlem{7.1}
\def \BIDUALone{27}
\def \PERPthm{7.2}
\def \PERPone{28}
\def \BBcor{7.3}
\def \APPsec{8}
\def \FACTORlem{8.1}
\def \ARENS{1}
\def \BB{2}
\def \BBWYLIN{3}
\def \BBWYFP{4}
\def \BBWYBB{5}
\def \BRON{6}
\def \GMS{7}
\def \ASD{8}
\def \MLMAX{9}
\def \FENCHEL{10}
\def \RANGE{11}
\def \PANDM{12}
\def \HBM{13}
\def \SSDMON{14}
\def \SSDB{15}
\def \YAO{16}
\def \ZBOOK{17}
%
\magnification 1200
\headline{\ifnum\folio=1
{\hfil{\largebf Linear {\largemi L}--positive sets and their polar subspaces}\hfil}
\else\centerline{\rm {\bf Linear $L$--positive sets and their polar subspaces}}\fi}
\bigskip
\centerline{\bf S. Simons}
\bigskip
\centerline{\sl Department of Mathematics, University of California}
\centerline{\sl Santa Barbara, CA 93106-3080, U.S.A.}
\centerline{\sl simons@math.ucsb.edu}
\medskip
\centerline{\bf Abstract}
\medskip
\noindent
In this paper, we define a Banach SNL space to be a Banach space with a certain kind of linear map from it into its dual, and we develop the theory of linear $L$--positive subsets of Banach SNL spaces with Banach SNL dual spaces.   We use this theory to give simplified proofs of some recent results of Bauschke, Borwein, Wang and Yao, and also of the classical Brezis--Browder theorem.
\defSection \INTsec
\medbreak
\centerline{\bf \INTsec.\quad Introduction}
\medskip
\noindent
This paper was motivated by a recent result of Bauschke, Borwein, Wang and Yao (which generalizes the Brezis--Browder theorem for reflexive spaces) that \slant if $A$ is a norm--closed linear monotone subspace of the product of a real Banach space with its dual space then the adjoint, $A^*$, of $A$ is monotone if, and only if, $A$ is maximally monotone of type (NI)\endslant.   This result appear in Corollary \MONcor.   We will show in Theorem \INFthm\ how to obtain relatively simple proofs of generalizations of this results in the context of Banach SNL spaces using the basic tools of convex analysis (that is to say Rockafellar's formula for the subdifferential of a sum, and the Br\o ndsted--Rockafellar theorem).
\par    
We introduce Banach SNL spaces in Section \SNLsec, and in Section \MONNIsec, we give details of how this concept can be applied to monotone subsets of the product of a real Banach space with its dual space.   In particular, we show in Lemma \RECASTlem\ how the definition of the concept of \slant maximal monotone subset of type (NI)\endslant\ can be recast in a very simple way in the notation of Banach SNL spaces. 
\par
In Section \CONVsec, we lay out for the benefit of the reader the basic tools of convex analysis (mentioned above) that we will use.
\par
The central section of this paper, Section \LINsec, is about linear $L$--positive sets and their polar subspaces.   If $C$ is a convex  $L$--positive subset of $B$, we introduce a function $q^C$; we show in Lemma \QCSTARlem\ that this function is convex, and we characterize its subdifferential in Lemma \QClem.  Theorem \INFthm(a) is the main result in this paper, and Theorem \INFthm(b--d) are all fairly simple consequences of Theorem \INFthm(a).  
\par
In Section \LINMONsec, we show how to deduce the results on linear monotone subspaces of the product of a real Banach space with its dual that we mentioned above and, in Section \SURBBsec, we introduce a surjectivity hypothesis that enables us to deduce the Brezis--Browder theorem for reflexive spaces.
\par
In Section \APPsec, we describe an isometry result on Banach SNL spaces with Banach SNL duals that was not needed in the previous sections.  
\par
While it is true that the introduction of SNL spaces to solve certain problems on monotone sets does increase the level of abstraction, it is also true that the resulting simplification of the proofs makes it easier to see the essential structure behind them.  In particular, the inequality (\PTtwo) \big(which appears again in (\INFCzero)\big) turns out to be very important.
\par
The author would like to thank Liangjin Yao for sending him the preprints \cite\BBWYLIN\endcite, \cite\BBWYFP\endcite, \cite\BBWYBB\endcite\ and \cite\YAO\endcite, and Maicon Marques Alves for sending him the preprint \cite\ASD\endcite.
\defSection\SNLsec
\medbreak
\centerline{\bf \SNLsec.\quad Banach SNL spaces and Banach SNL duals}
\medskip
\noindent
We start off by introducing some Banach space notation.
\defDefinition \PCdef
\medbreak
\noindent
{\bf Definition \PCdef.}\enspace If $X$ is a nonzero real Banach space and $f\colon\ X \to \rbar$, we say that $f$ is \slant proper\endslant\ if there exists $x \in X$ such that $f(x) \in \r$.   We write $X^*$ for the dual space of $X$ \big(with the pairing $\bra\cdot\cdot\colon X \times X^* \to \r$\big) and $X\dbs$ for the bidual of $X$ \big(with the pairing $\bra\cdot\cdot\colon X^* \times X\dbs \to \r$\big).   If $x \in X$, we write $\wh x$ for the canonical image of $x$ in $X\dbs$, that is to say\quad $x \in X\ \and\ x^* \in X^* \lr\bra{x^*}{\wh x} = \bra{x}{x^*}$.  
\medskip
We now introduce the concept of Banach SNL space.   Many of the results in this section appear (implicitly) in \cite\PANDM\endcite, \cite\HBM\endcite, \cite\SSDMON\endcite\ and \cite\SSDB\endcite.
\defDefinition \NORMdef
\medskip
\noindent
{\bf Definition \NORMdef.}\enspace Let $B$ be a nonzero real Banach space.   A \slant SNL map on $B$\endslant\  (``SNL'' stands for ``symmetric nonexpansive linear''), is a linear map $L\colon\ B \to B^*$ such that
$$\|L\| \le 1\quad\and,\quad\all\ b,c\in B,\quad \bra{b}{Lc} = \bra{c}{Lb}.\meqno\IOTAone$$
A \slant Banach SNL space\enspace $(B,L)$\endslant\ is a nonzero real Banach space $B$ together with a SNL map $L\colon\ B \to B^*$.   We define the function $q_L\colon\ B \to \r$ by\quad $q_L(b) := \half\bra{b}{Lb}$\quad($b \in B$)\quad   (``$q$'' stands for ``quadratic'').   Clearly
$$q_L + \half\|\cdot\|^2 \ge 0\ \on\ B.\meqno\Pone$$
It also follows from (\IOTAone) that,\quad for all $d,e \in B$,\quad$|q_L(d) - q_L(e)| = \half|\bra{d}{Ld} - \bra{e}{Le}| =\break\half\big|\Bra{d - e}{L(d + e)}\big| \le \half\|d - e\|\|d + e\|$, and so
$$q_L\ \hbox{is continuous}.\meqno\QCONTone$$
Now let $(B,L)$ be a Banach SNL space and $A \subset B$.   We say that $A$ is \slant$L$--positive\endslant\ if $A \ne \emptyset$ and\quad $b,c \in A \lr q_L(b - c) \ge 0$.\quad We say that $A$ is \slant$L$--negative\endslant\ if $A \ne \emptyset$ and \quad $b,c \in A \lr q_L(b - c) \le 0$.\quad   We say that $A$ is \slant maximally $L$--positive\endslant\ if $A$ is $L$--positive and $A$ is not properly contained in any other $L$--positive set.   In this case,
$$d \in B \lr \infn_{a \in A}q_L(d - a) \le 0.\meqno\QMAXone$$
(If $d \in B \setminus A$ then the maximality gives us $a \in A$ such that $q_L(d - a) < 0$, while if $d \in A$ then $q_L(d - d) = 0$.)   Similarly, we say that $A$ is \slant maximally $L$--negative\endslant\ if $A$ is $L$--negative and $A$ is not properly contained in any other $L$--negative set.
\par
Now let $(B,L)$ be a Banach SNL space and $\big(B^*,\Lt\big)$ also be a Banach SNL space.   We say that $\big(B^*,\Lt\big)$ is a \slant Banach SNL dual\endslant\ of $(B,L)$ if,
$$\all\ b \in B,\qquad \Lt(Lb) = \wh b.\meqno\DUALone$$
Then (\DUALone) and (\IOTAone) imply that,
$$\all\ b \in B \ \and\ b^* \in B^*,\quad\qLt(b^* + Lb) = \half\Bra{b^* + Lb}{\Lt b^* + \wh b} = \qLt(b^*) + \bra{b}{b^*} + q_L(b).\meqno\ITone$$
In particular,
$$\qLt \circ L = q_L.\meqno\QTone$$
By analogy with (\Pone), we have
$$\qLt + \half\|\cdot\|^2 \ge 0\ \on\ B^*.\meqno\DUALthree$$
If $(B,L)$ is a Banach SNL space then so also is $(B,-L)$.   Further, if $\big(B^*,\Lt\big)$ is a Banach SNL dual of $(B,L)$ then, as the reader can easily verify, \slant$\big(B^*,-\Lt\big)$ is a Banach SNL dual of $(B,-L)$\endslant.
\medskip
There are many examples of Banach SNL spaces and their associated $L$--positive sets.   The following are derived from \cite\SSDMON, Examples 2.3, 2.5, Remarks 6.3, 6.7, 6.8, pp.\ 230--231, 244--246\endcite.   More examples can be found in \cite\GMS\endcite. 
\defExample \Hex 
\medbreak
\noindent
{\bf Example \Hex.}\enspace Let $B$ be a Hilbert space with inner product $(b,c) \mapsto \bra{b}{c}$ and $L\colon B \to B$ be a nonexpansive self--adjoint linear operator.   Then $(B,L)$ is a Banach SNL space.   Here are three special cases of this example:
\smallbreak
(a)\enspace If, for all $b \in B$, $Lb = b$ then every subset of $B$ is $L$--positive.
\smallbreak
(b)\enspace If, for all $b \in B$, $Lb = -b$ then the $L$--positive sets are the singletons.
\smallbreak
(c)\enspace If $B = \r^3$ and $L(b_1,b_2,b_3) = (b_2,b_1,b_3)$ and $M$ is any nonempty monotone subset of $\r \times \r$ (in the obvious sense) then $M \times \r$ is an $L$--positive subset of $B$.   The set $\r(1,-1,2)$ is an $L$--positive subset of $B$ which is not contained in a set $M \times \r$ for any monotone subset of $\r \times \r$.   The helix $\big\{(\cos\theta,\sin\theta,\theta)\colon \theta \in \r\big\}$ is an $L$--positive subset of $B$, but if $0 < \lambda < 1$ then the helix $\big\{(\cos\theta,\sin\theta,\lambda\theta)\colon \theta \in \r\big\}$ is not.
\smallbreak
(d)\enspace If $B = \r^3$ and $L(b_1,b_2,b_3) = (b_2,b_3,b_1)$ then $(B,L)$ cannot be a Banach SNL space, since (\IOTAone) fails. 
\medskip
In this paper, we are particularly interested in \cite\SSDMON, Example 6.5, p.\ 245\endcite, which we will explain in Example \EEex\ below.
\defSection \MONNIsec
\medbreak
\centerline{\bf \MONNIsec\quad  Connection with monotone sets:  type (NI)}
\medskip
\noindent
We suppose in this section that $E$ is a nonzero Banach space.
\defExample \EEex 
\medbreak
\noindent
{\bf Example \EEex.}\enspace Let $B := E \times E^*$ and, for all $(x,x^*) \in B$, $\|(x,x^*)\| := \sqrt{\|x\|^2 + \|x^*\|^2}$.   We represent $B^*$ by $E^* \times E\dbs$, under the pairing $\Bra{(x,x^*)}{(y^*,y\dbs)} := \bra{x}{y^*} + \bra{x^*}{y\dbs}$, and define $L\colon\ B \to B^*$ by $L(x,x^*) := (x^*,\wh{x})$.   Then $(B,L)$ is a Banach SNL space and, for all $(x,x^*) \in B$, $q_L(x,x^*) = \half\big[\bra{x}{x^*} + \bra{x}{x^*}\big] = \bra{x}{x^*}$.   If $(x,x^*), (y,y^*) \in B$ then we have $q_L\big((x,x^*) - (y,y^*)\big) = q_L(x - y,x^* - y^*) = \bra{x - y}{x^* - y^*}$.   Thus if $A \subset B$ then $A$ is $L$--positive exactly when $A$ is a nonempty monotone subset of $B$ in the usual sense, and $A$ is maximally $L$--positive exactly when $A$ is a maximally monotone subset of $B$ in the usual sense.   We point out that any finite dimensional Banach SNL space of the form described here must have \slant even\endslant\ dimension, and that there are many Banach SNL spaces  of finite odd dimension with Banach SNL duals.   See \cite\SSDMON, Remark 6.7, p.\ 246\endcite.
\par
As usual, the dual norm on $B^* = E^* \times E\dbs$ is given by  $\|(y^*,y\dbs)\| := \sqrt{\|y^*\|^2 + \|y\dbs\|^2}$. By analogy with the analysis above, we define $\Lt\colon\ B^* \to B\dbs$ by $\Lt(y^*,y\dbs) = \big(y\dbs,\wh{y^*}\big)$. Then $\big(B^*,\Lt\big)$ is a Banach SNL space, and, for all $(y^*,y\dbs) \in B^*$, $\qLt(y^*,y\dbs) = \bra{y^*}{y\dbs}$.    Now let $(x,x^*) \in B$.   Then $\Lt\big(L(x,x^*)\big) = \Lt(x^*,\wh{x}) = \big(\wh{x},\wh{x^*}\big) = \wh{(x,x^*)}$.   Thus (\DUALone) is satisfied, and so $\big(B^*,\Lt\big)$ is a Banach SNL dual of $(B,L)$.
\par
Finally, let $d^* =(y^*,y\dbs) \in B^*$ and $\eps > 0$.   From the definition of $\|y\dbs\|$, there exists $z^* \in E^*$ such that $\|z^*\| \le \|y\dbs\|$ and $\bra{z^*}{y\dbs} \ge \|y\dbs\|^2 - \eps$.   Let $d = (0,y^* + z^*) \in B$.   But then $d^* - Ld = (y^*,y\dbs) - (y^* + z^*,0) = (-z^*,y\dbs) \in B^*$.   Since
$$\qLt(-z^*,y\dbs) + \half\big\|(-z^*,y\dbs)\big\|^2 =   - \bra{z^*}{y\dbs} + \half\big(\|z^*\|^2 + \|y\dbs\|^2\big) \le  - \bra{z^*}{y\dbs} + \|y\dbs\|^2 \le \eps,$$
we have shown that
$$\all\ d^* \in B^*,\quad \infn_{d \in B}\big[\qLt(d^* - Ld) + \half\|d^* - Ld\big\|^2] \le 0.\meqno\PTtwo$$\par
\defRemark \PTrem
\par
\noindent
{\bf Remark \PTrem.}\enspace In \cite\SSDMON, Definition 6.4, p.\ 245\endcite, we used the terminology ``$L(B)$ is \slant $\pt$--dense in\endslant\ $B^*$'' to describe the inequality (\PTtwo).
\medskip
Definition \NIdef\ first appeared in \cite\RANGE, Definition 10, p.\ 183\endcite.   It was thought at first that this was a weak definition, but it was proved by Marques Alves and Svaiter in \cite\ASD, Theorem 4.4, pp.\ 1084--1086\endcite, that \slant if $A$ is a maximally monotone subset of $E \times E^*$ of type (NI) then $A$ is of ``type (D)''\endslant\ (the opposite implication is obviously true).   This result was extended in \cite\SSDMON, Theorem 9.9(a), pp.\ 254--255\endcite, where it was proved that \slant if $A$ is a maximally monotone subset of $E \times E^*$ of type (NI) then $A$ is of ``dense type'' and ``type (ED)''\endslant.   This has a number of consequences, which are detailed in \cite\SSDMON, Theorem 9.9(b--f)\endcite\ and \cite\SSDMON, Theorem 9.10\endcite.   Finally, it was proved by Bauschke, Borwein, Wang and Yao in \cite\BBWYFP, Theorem 3.1\endcite\ that \slant if $A$ is a maximally monotone subset of $E \times E^*$ then $A$ is of type (NI) if, and only if, $A$ is of ``Fitzpatrick--Phelps'' type\endslant.
\medskip
What distinguishes the definition of ``type (NI)'' from the definitions of the other classes of maximally monotone sets mentioned above is that it can be recast easily in the language of Banach SNL spaces, as we will see in Lemma \RECASTlem\ below.  
\defDefinition \NIdef
\medbreak
\noindent
{\bf Definition \NIdef.}\enspace Let $A \subset E \times E^*$.   We say that $A$ is \slant maximally monotone of type (NI)\endslant\ if $A$ is maximally monotone and,
$$\all\ (y^*,y\dbs) \in E^* \times E\dbs,\quad \infn_{(x,x^*) \in A}\bra{y^* - x^*}{y\dbs - \wh x} \le 0.$$
\par
Here is the promised reformulation of Definition \NIdef\ in the notation of Banach SNL spaces.   We use the conventions introduced in Example \EEex.
\defLemma \RECASTlem
\medbreak
\noindent
{\bf Lemma \RECASTlem}\enspace\slant Let $A \subset B := E \times E^*$.   Then $A$ is maximally monotone of type (NI) if, and only if, $A$ is maximally $L$--positive and, for all $d^* \in B^*$, $\infn_{a \in A}\qLt(d^* - La) \le 0$.\endslant
\Proof This is immediate from the formulae for $\qLt$ and $L$ given in Example \EEex.\qed
\defSection \CONVsec
\medbreak
\centerline{\bf \CONVsec\quad  Standard results from convex analysis}
\medskip
\noindent
For the rest of this paper, we suppose that $(B,L)$ is a Banach SNL space and that $\big(B^*,\Lt\big)$ is a Banach SNL dual of $(B,L)$.   We will use the following standard notation and results from convex analysis.   If $f\colon\ B \to \rbar$ is proper and convex and $b^* \in B^*$ then the \slant Fenchel conjugate\endslant, $f^*$, of $f$ is defined by\
$$f^*(b^*) := \supn_B\big[b^* - f\big].$$
If, further, $b \in B$ then the \slant subdifferential of $f$ at $b$\endslant, $\partial f(b)$, is the subset of $B^*$ defined by
$$b^* \in \partial f(b) \iff f(b) + f^*(b^*) = \bra{b}{b^*}.$$
We recall \big(see Rockafellar, \cite\FENCHEL, Theorem 3(b), pp.\ 85--86\endcite, Z\u alinescu, \cite\ZBOOK, Theorem 2.8.7(iii), p.\ 127\endcite, or \cite\HBM, Theorem 18.1, p.\ 74\endcite\big) the following special case of:
\medbreak
\noindent
{\bf Rockafellar's formula for the subdifferential of a sum.}\enspace\slant Let $f\colon\ B \to \rbar$ be proper and convex.   Then
$$\partial(f + \half\|\cdot\|^2) = \partial f + J,\meqno\Jone$$
where $J = \partial\big(\half\|\cdot\|^2\big)$ is characterized by
$$c^* \in Jc \iff \|c\|^2 = \|c^*\|^2 = \bra{c}{c^*}.\meqno\Jtwo$$\endslant
$J$ is called the ``duality map''.   Finally, we recall \big(see Br\o ndsted--Rockafellar, \cite\BRON, p.\ 608\endcite, Z\u alinescu, \cite\ZBOOK, Theorem 3.1.2, p. 161\endcite, or \cite\HBM, Theorem 18.6, p.\ 76\endcite\big) the following special case of the:
\medbreak
\noindent
{\bf Br\o ndsted--Rockafellar theorem.}\enspace \slant Let $f\colon\ B \to \rbar$ be proper convex and lower semicontinuous, $\eta > 0$, $b_0 \in B$, $b_0^* \in B^*$ and\quad $f(b_0) + f^*(b_0^*) \le \bra{b_0}{b_0^*} + \eta^2$.\quad Then there exist $b_1 \in B$ and $b_1^* \in \partial f(b_1)$ such that $\|b_1 - b_0\| \le \eta$ and $\|b_1^* - b_0^*\| \le \eta$.\endslant
\defSection \LINsec
\medbreak
\centerline{\bf \LINsec\quad  Linear $L$--positive sets}
\medskip
\noindent
In this section, we suppose that $A$ is a linear $L$--positive subspace of $B$. If $C$ is a convex $L$--positive subset of $B$, we define the functions\quad $q^C,p^C\colon\ B \to \,\rbar$\quad by
$$q^C(b) = \cases{q_L(b)&(if $b \in C$);\cr \infty&(if $b \in B \setminus C$),}\quand p^C(b):= \cases{q_L(b)+ \half\|b\|^2&(if $b \in C$);\cr \infty&(if $b \in B \setminus C$).}$$
Now fix $c_0 \in C$.   Then, for all $b \in C$, $q_L(b) + \|b\|\|c_0\| + q_L(c_0) \ge q_L(b) - \bra{b}{Lc_0} + q_L(c_0) = q_L(b - c_0) \ge 0$ thus, for all $b \in B$, $p^C(b) = q^C(b) + \half\|b\|^2 \ge - \|b\|\|c_0\| - q_L(c_0) + \half\|b\|^2$.   Consequently
$$p^C(b) \to \infty\quad \hbox{as}\quad \|b\| \to \infty.\meqno\COERCIVITY$$\par
\defLemma \QCSTARlem
\noindent
{\bf Lemma \QCSTARlem}\enspace\slant Let $A$ be a linear $L$--positive subspace of $B$,  $d \in B$, $C := A - d$ and $d^* \in B^*$.  Then $q^C$ is proper and convex.\endslant
\Proof $q^C$ is obviously proper.   Suppose that $b,c \in C$ and $\lambda \in \,]0,1[\,$.   Then
$$\lambda q_L(b) + (1 - \lambda)q_L(c) - q\big(\lambda b + (1 - \lambda)c\big) = \lambda(1 - \lambda)q_L(b - c) \ge 0.$$
This implies the convexity of $q^C$.  \big(See \cite\HBM, Lemma 19.7, pp.\ 80--81\endcite.\big)\qed
\medbreak
We write $A^0$ for the linear subspace  $\big\{b^* \in B^*\colon\ \bra{A}{b} = \{0\}\big\}$ of $B^*$.   $A^0$ is the ``polar subspace of $A$''.   The significance of $A^0$ lies in the following lemma:  
\defLemma \QClem
\medskip
\noindent
{\bf Lemma \QClem}\enspace\slant Let $A$ be a linear $L$--positive subspace of $B$.  $d \in B$ and $C := A - d$.  Then
$$\partial q^C(b) = \cases{Lb + A^0&(if $b \in C$);\cr \emptyset&(if $b \in B \setminus C$).}$$\endslant
\Proo Since it is obvious that\quad $\partial q^C(b) = \emptyset$ if $b \in B \setminus C$,\quad it only remains to show that\quad $\partial q^C(b) = Lb + A^0$ if $b \in C$.\quad So suppose that $b \in C$.   Then, since $c - b$ runs through $A$ as $c$ runs through $C$ and, from (\IOTAone),\quad $q_L(b) - q_L(a + b) = -\Bra{a}{Lb} - q_L(a)$,
$$\eqalign{b^* \in \partial q^C(b)
&\iff q^C(b) + \supn_{c \in C}\big[\bra{c}{b^*} - q_L(c)\big] \le \bra{b}{b^*}\cr
&\iff \all\ c \in C,\ \bra{c - b}{b^*} + q_L(b) - q_L(c) \le 0\cr
&\iff \all\ a \in A,\ \bra{a}{b^*} + q_L(b) - q_L(a + b) \le 0\cr
&\iff \all\ a \in A,\ \bra{a}{b^* - Lb} \le q_L(a).}$$
Since $q_L(a) \ge 0$, this is trivially satisfied if $b^* \in Lb + A^0$.   On the other hand, if $b^* \in \partial q^C(b)$ then it follows from the above that, for all $a \in A$ and $\lambda \in \r$, 
$$\lambda \bra{a}{b^* - Lb} = \bra{\lambda a}{b^* - Lb} \le q_L(\lambda a) = \lambda^2q_L(a),$$
and, from standard quadratic arguments, $b^* - Lb \in A^0$, that is to say $b^* \in Lb + A^0$.\qed
\medskip
Theorem \INFthm\ was suggested by the analysis in Bauschke, Borwein, Wang and Yao\break \cite\BBWYLIN, Theorem 3.1(iii)$\lr$(ii)\endcite\ and \cite\BBWYBB, Proposition 3.1\endcite.   In contrast to the analysis in these two papers, the only results on convex functions that we use are Rockafellar's formula for the conjugate of a sum, and the Br\o ndsted--Rockafellar theorem, as outlined in Section \CONVsec\ above.   
\defTheorem \INFthm
\medbreak
\noindent
{\bf Theorem \INFthm.}\enspace\slant Suppose that $A$ be a norm--closed linear $L$--positive subspace of $B$ and $A^0$ is $\Lt$--negative.   Then:
\medskip
\noindent
{\rm(a)}\enspace Let $d \in B$ and $C := A - d$.   Then $\infn_Cp^C \le 0$. 
\medskip
\noindent
{\rm(b)}\enspace Let $d \in B$, $C := A - d$ and $b^* \in B^*$.   Then $-\big({q^C}\big)^*(b^*) \le \half\|b^*\|^2$. 
\medskip
\noindent
{\rm(c)}\enspace Let $d\in B$ and $d^* \in B^*$.   Then $\infn_{a \in A}\qLt(d^* - La) \le \qLt(d^* - Ld) + \half\|d^* - Ld\|^2$. 
\medskip
\noindent
{\rm(d)}\enspace $A$ is maximally $L$--positive.
\endslant
\Proof(a)\enspace It is clear from (\Pone) and (\COERCIVITY) that there exists $M > 0$ such that, writing $C_M := \{b\colon\ b \in C,\ \|b\| \le M\big\}$, $\infn_{C_M}p^C = \infn_Cp^C \in [\,0,\infty[\,$.   Now let $\eps > 0$ be arbitrary. Choose $\eta > 0$ so that $2M\eta + 3\eta^2 \le \eps$, and $b \in C_M$ such that $p^C(b) \le \infn_{C_M}p^C + \eta^2 = \infn_Cp^C + \eta^2 = \infn_Bp^C + \eta^2$.   We have
$$\|b\| \le M \quand p^C(b) + \big({p^C}\big)^*(0) - \bra{b}{0} = p^C(b) - \infn_Bp^C \le \eta^2 .\meqno\INFone$$
From (\QCONTone), $q_L$ is continuous.   Since $C$ is norm--closed in $B$, $p^C$ is lower semicontinous, and so the Br\o ndsted--Rockafellar theorem provides us with $c \in B$ and $b^* \in \partial p^C(c)$ such that $\|c - b\| \le \eta$ and $\|b^*\| \le \eta$.   From Rockafellar's formula for the subdifferential of a sum, (\Jone), $\partial p^C = \partial q^C + J$.  Combining this with (\INFone) and Lemma \QClem, we have
$$\|c\| \le M + \eta,\ \|b^*\| \le \eta,\ c \in C\ \and\ \exs\ c^* \in Jc\ \st\ b^* \in Lc + A^0 + c^*.\meqno\INFtwo$$
From (\IOTAone), (\Jtwo) and (\INFtwo), we have
$$\|b^*\| \le \eta \quand \|Lc + c^*\| \le \|Lc\| + \|c^*\| \le 2\|c\| \le 2(M + \eta)\meqno\INFthree$$
and, from (\Jtwo) and (\ITone), we have
$$q_L(c) + \|c\|^2 + \qLt(c^*) = q_L(c) + \bra{c}{c^*} + \qLt(c^*) = \qLt\big(Lc + c^*\big).\meqno\INFfour$$
(\INFtwo) implies that $b^* - Lc - c^* \in A^0$. From (\Jtwo), (\DUALthree), (\INFfour), the $\Lt$--negativity of $A^0$ and (\INFthree),
$$\eqalign{p^C(c) &= q_L(c) + \half\|c\|^2 = q_L(c) + \|c\|^2 - \half\|c^*\|^2 \le q_L(c) + \|c\|^2 + \qLt(c^*) = \qLt\big(Lc + c^*\big)\cr
&= \qLt\big(b^* - Lc - c^*\big) + \Bra{Lc + c^*}{\Lt b^*} - \qLt(b^*)
\le \Bra{Lc + c^*}{\Lt b^*} + \half\|b^*\|^2 \cr
&\le \|b^*\|\|Lc + c^*\| + \half\|b^*\|^2 \le 2(M + \eta)\eta + \half\eta^2 < \eps.}$$
Since $\eps > 0$ is arbitrary, this completes the proof of (a).
\medskip
(b)\enspace From the definition of Fenchel conjugate (or the Fenchel--Young inequality), for all $c \in C$,\quad $-\big({q^C}\big)^*(b^*) \le q^C(c) - \bra{c}{b^*} \le q^C(c) + \|c\|\|b^*\| \le q^C(c) + \half\|c\|^2 + \half\|b^*\|^2 = p^C(c) + \half\|b^*\|^2$,\quad and the result now follows from (a).
\medskip
(c)\enspace Let $C := A - d$.   Using (\ITone), and the definition of $q^C$, we have 
$$\eqalign{\infn_{a \in A}\qLt(d^* - La)
&= \infn_{a \in A}\big[\qLt(d^*) - \bra{a}{d^*} + q_L(a)\big]\cr
&= \infn_{c \in C}\big[\qLt(d^*) - \bra{d + c}{d^*} + q_L(d + c)\big]\cr
&= \infn_{c \in C}\big[\qLt(d^*) - \bra{d}{d^*} - \bra{c}{d^*} + q_L(d) + \bra{c}{Ld} + q_L(c)\big]\cr
&= \infn_{c \in C}\big[\qLt(d^* - Ld) - \bra{c}{d^* - Ld} + q_L(c)\big]\cr
&= \qLt(d^* - Ld) - \supn_{c \in C}\big[\bra{c}{d^* - Ld} - q^C(c)\big]\cr
&= \qLt(d^* - Ld) - \big({q^C}\big)^*(d^* - Ld),
}$$
and the result follows from (b) with $b^* := d^* - Lc$.
\medskip 
(d)\enspace Let $d \in B$ and $A \cup \{d\}$ be $L$--positive.   Then, for all $a \in A$, $q_L(a - d) \ge 0$.   If now $C := A - d$ then $q^C \ge 0$ on $C$, and so $p^C = q^C + \half\|\cdot\|^2 \ge \half\|\cdot\|^2$ on $C$.   Thus (a) implies that $\inf_C\half\|\cdot\|^2 \le 0$, from which $0 \in \overline{C}$, which implies in turn that $d \in \overline{A}$.   Since $A$ is closed, $d \in A$.   This completes the proof that $A$ is maximally $L$--positive.\qed
\par
\defCorollary \INFcor
\medbreak
\noindent
{\bf Corollary \INFcor.}\enspace\slant Let $A$ be a norm--closed linear $L$--positive subspace of $B$ and
$$\all\ d^* \in B^*,\quad \infn_{d \in B}\big[\qLt(d^* - Ld) + \half\|d^* - Ld\|^2\big] \le 0.\meqno\INFCzero$$
Then\endslant\ (\INFCone)$\iff$(\INFCtwo)$\iff$(\INFCthree)$\lr$(\INFCfour).\slant
$$A^0\ \hbox{is}\ \Lt\hbox{--negative.}\meqno\INFCone$$
$$\All\ d^* \in B^*,\ \infn_{a \in A}\qLt(d^* - La) \le 0.\meqno\INFCtwo$$
$$\All\ d^* \in A^0,\ \infn_{a \in A}\qLt(d^* - La) \le 0.\meqno\INFCthree$$
$$A\ \hbox{is maximally}\ L\hbox{--positive.}\meqno\INFCfour$$\endslant
\Proo It is immediate from Theorem \INFthm(c) and (\INFCzero) that (\INFCone)$\lr$(\INFCtwo), and it is obvious that (\INFCtwo)$\lr$(\INFCthree).  If (\INFCthree) is true then, using (\ITone), for all $d^* \in A^0$,
$$\qLt(d^*) + \infn_{a \in A}q_L(a) = \infn_{a \in A}\big[\qLt(d^*) - \bra{a}{d^*} + q_L(a)\big] = \infn_{a \in A}\qLt(d^* - La) \le 0.$$
Since $\infn_{a \in A}q_L(a) \ge 0$, it follows that $\qLt(d^*) \le 0$, and the linearity of $A^0$ gives (\INFCone).   Finally, it is immediate from Theorem \INFthm(d) that (\INFCone)$\lr$(\INFCfour).\qed
\medskip
\defCorollary \NEGcor
\medbreak
\noindent
{\bf Corollary \NEGcor.}\enspace\slant If $A$ is a norm--closed linear $L$--negative subspace of $B$, and $A^0$ is $\Lt$--positive then $A$ is maximally $L$--negative.\endslant
\Proof This result follows by applying Theorem \INFthm(d) to $(B,-L)$, since\break $L$--negativity is equivalent to $(-L)$--positivity and, $\Lt$--positivity is equivalent to\break $\big(-\Lt\big)$--negativity.\qed
\medskip
\defSection \LINMONsec
\medbreak
\centerline{\bf \LINMONsec.\quad Linear monotone subspaces}
\medskip
\noindent
We suppose in this section that $E$ is a nonzero Banach space, and we use the notation of Example \EEex.   We define $\rho\colon\ E^* \times E\dbs \to E^* \times E\dbs$ by $\rho(x^*,x\dbs) := (x^*,-x\dbs)$.   Let $\emptyset \ne D \subset E^* \times E\dbs$.   Since $\qLt \circ \rho = -\qLt$,
$$D\ \hbox{is}\ \Lt\hbox{--negative} \iff \rho(D)\ \hbox{is}\ \Lt\hbox{--positive}.\meqno\RHOone$$
If $A$ is linear subspace of $E \times E^*$ then the \slant adjoint subspace\endslant, $A^*$, of $E^* \times E\dbs$, is defined by:\quad  $(y^*,y\dbs) \in A^* \iff \hbox{for all}\ (a,a^*) \in A,\ \bra{a}{y^*} = \bra{a^*}{y\dbs}$.\quad This definition goes back at least to Arens in \cite\ARENS\endcite.   Clearly
$$A^* = \rho\big(A^0\big).\meqno\ADJPOLone$$
\par
The implication (\MONone)$\lr$(\MONfour) below was shown by Bauschke, Borwein, Wang and Yao in a two--stage process as follows:   It was first shown in \cite\BBWYLIN, Theorem 3.1(iii)$\lr$(ii)\endcite\ that if $A$ is a maximally monotone linear subspace of $E \times E^*$ and (\MONone) is satisfied then $A$ is of type (NI), and it was subsequently shown in \cite\BBWYBB, Proposition 3.1\endcite\ that if $A$ is a norm--closed monotone linear subspace of $E \times E^*$ and (\MONone) is satisfied then $A$ is maximal.  
\defCorollary \MONcor
\medbreak
\noindent
{\bf Corollary \MONcor.}\enspace\slant Let $A$ be a norm--closed linear monotone subspace of $E \times E^*$.  Then the four conditions below are equivalent.
$$A^*\ \hbox{is monotone.}\meqno\MONone$$
$$\All\ (y^*,y\dbs) \in E^* \times E\dbs,\quad \infn_{(x,x^*) \in A}\bra{y^* - x^*}{y\dbs - \wh x} \le 0.$$
$$\All\ (y^*,y\dbs) \in A^0,\quad \infn_{(x,x^*) \in A}\bra{y^* - x^*}{y\dbs - \wh x} \le 0.$$
$$A\ \hbox{is maximally monotone of type (NI).}\meqno\MONfour$$\endslant
\Proo This is immediate from (\ADJPOLone), (\RHOone), (\PTtwo), Corollary \INFcor, and Lemma \RECASTlem.\qed
\defSection \SURBBsec
\medbreak
\centerline{\bf \SURBBsec.\quad The surjectivity hypothesis and the Brezis--Browder theorem}
\medskip
\noindent
In this section, we will consider the implications of the \slant surjectivity hypothesis\endslant:
$$L(B) = B^*.\meqno\SURone$$\par
\defLemma \SURlem
\par
\noindent
{\bf Lemma \SURlem.}\enspace\slant Let $(B,L)$ be a Banach SNL space, $\big(B^*,\Lt\big)$ be a Banach SNL dual of $(B,L)$ and  {\rm(\SURone)} be satisfied.   Then (writing $\cdot^{\rm T}$ for the adjoint of $\cdot$ in the usual sense):
\par
\noindent
{\rm(a)}\enspace $\big(B\dbs,\Ltt\big)$ is a Banach SNL dual of $(B^*,\Lt)$.
\par
\noindent
{\rm(b)}\enspace $B$ is reflexive.\endslant
\Proof Clearly, $\Ltt\colon\ B\dbs \to B^{***}$ is linear and $\big\|\Ltt\big\| =  \big\|L^{\rm T}\big\| = \|L\| \le 1$.   Now let $b\dbs, c\dbs \in B\dbs$.   Then $L^{\rm T}b\dbs, L^{\rm T}c\dbs \in B^*$, and (\SURone) provides us with $b,c \in B$ such that\quad $Lb = L^{\rm T}b\dbs$\quad and \quad $Lc = L^{\rm T}c\dbs$.\quad   Thus, from (\IOTAone),
$$\eqalign{&\Bra{b\dbs}{\Ltt c\dbs} = \Bra{L^{\rm T}b\dbs}{c\dbs} = \bra{Lb}{c\dbs}  = \Bra{b}{L^{\rm T}c\dbs} = \bra{b}{Lc}\cr
&\qquad= \bra{c}{Lb} = \Bra{c}{L^{\rm T}b\dbs} = \bra{Lc}{b\dbs} = \Bra{L^{\rm T}c\dbs}{b\dbs} = \Bra{c\dbs}{\Ltt b\dbs}.}$$
Thus the analog of (\IOTAone) is satisfied, and so $\big(B\dbs,\Ltt\big)$ is a Banach SNL space.   We next show that
$$\all\ c^* \in B^*,\qquad \Ltt\big(\Lt c^*\big) = \wh{c^*} \in  B^{***}.\meqno\BIDUALone$$
To this end, let $c^* \in B^*$, and $b\dbs$ be an arbitrary element of $B\dbs$.  (\SURone) provides us with $c \in B$ such that\quad $Lc = c^*$.\quad   Then, from (\DUALone),\quad $\Lt c^* =  \Lt(Lc) = \wh c$,\quad and so
$$\eqalign{\Bra{b\dbs}{\Ltt\big(\Lt c^*\big)}
&= \Bra{L^{\rm T}b\dbs}{\Lt c^*} = \Bra{L^{\rm T}b\dbs}{\wh c}\cr
&= \Bra{c}{L^{\rm T}b\dbs} = \Bra{Lc}{b\dbs} = \bra{c^*}{b\dbs}.}$$ 
Thus (\BIDUALone) holds.   Since (\BIDUALone) is the analog of (\DUALone), this completes the proof of (a).
\medbreak  
(b)\enspace Let $b\dbs \in B\dbs$.   Since\quad $L^{\rm T}b\dbs \in B^*$,\quad (\SURone) provides us with $b \in B$ such that\quad $Lb = L^{\rm T}b\dbs$.\quad  Now let $c^*$ be an arbitrary element of $B^*$.   From (\SURone) again, there exists $c \in B$ such that\quad $Lc = c^*$.\quad   Then, from (\IOTAone),
$$\bra{b}{c^*} = \Bra{b}{Lc} = \Bra{c}{Lb} = \Bra{c}{L^{\rm T}b\dbs} = \Bra{Lc}{b\dbs} = \bra{c^*}{b\dbs}.$$  
Thus $B$ is reflexive.   This completes the proof of (b).\qed
\medskip
Our next result is suggested by Yao, \cite\YAO, Theorem 2.4, p.\ 3\endcite.
\defTheorem \PERPthm
\medbreak
\noindent
{\bf Theorem \PERPthm.}\enspace\slant  Let the surjectivity hypothesis (\SURone) be satisfied, and $A$ be a norm--closed linear $L$--positive subspace of $B$.   Then:
\smallbreak\noindent
{\rm(a)}\enspace $A$ is maximally $L$--positive if, and only if, $A^0$ is $\Lt$--negative.
\smallbreak\noindent
{\rm(b)}\enspace $A$ is maximally $L$--positive if, and only if, $A^0$ is maximally $\Lt$--negative.\endslant
\Proof ``If'' in (a) is clear from Theorem \INFthm(d).   Conversely, let us suppose that $A$ is maximally $L$--positive.   Let $d^* \in B^*$.   From (\SURone), there exists $d \in B$ such that $Ld = d^*$.   Then, from (\QTone) and (\QMAXone),
$$\eqalign{\infn_{a \in A}\qLt(d^* - La)
&= \infn_{a \in A}\qLt(Ld - La)\cr
&= \infn_{a \in A}\qLt\big(L(d - a)\big) = \infn_{a \in A}q_L(d - a) \le 0,}$$   
and, since (\SURone) implies (\INFCzero), it follows from Corollary 5.4 that $A^0$ is $\Lt$--negative.   This completes the proof of (a).  ``If'' in (b) is immediate from ``If'' in (a).   Conversely, let us suppose that $A$ is maximally $L$--positive.   Now $A^0$ is norm--closed and, from (a), $\Lt$--negative, and it only remains to prove the $\Lt$--maximality.   We now show that
$$\all\ a\dbs \in (A^0)^0,\quad q_{L^{\rm TT}}(a\dbs) \ge 0.\meqno\PERPone$$   
To the end, let $a\dbs \in (A^0)^0$.   Since $A$ is norm--closed, standard functional analysis and the reflexivity of $B$ \big(Lemma \SURlem(b)\big) imply that there exists $a \in A$ such that $a\dbs = \wh a$.   From (\DUALone), Lemma \SURlem(a), and two applications of (\QTone) \big(one to $(B,L)$ and one to $\big(B^*,\Lt\big)$\big),
$$q_{L^{\rm TT}}(a\dbs) = q_{L^{\rm TT}}(\wh a) = q_{L^{\rm TT}}\big(\Lt(La)\big) = \qLt(La) = q_L(a) \ge 0.$$  
Thus (\PERPone) is true, from which it follows that $(A^0)^0$ is $L^{\rm TT}$--positive.   If we now apply Corollary \NEGcor, with $B$ replaced by $B^*$, and $A$ replaced by $A^0$, we see that $A^0$ is maximally $\Lt$--negative, which completes the proof of (b).\qed          
\medskip
Corollary \BBcor(a) below appears in Brezis--Browder \cite\BB, Theorem 2, pp.\ 32--33\endcite, and Corollary \BBcor(b) in Yao, \cite\YAO, Theorem 2.4, p.\ 3\endcite. 
\defCorollary \BBcor
\medbreak
\noindent
{\bf Corollary \BBcor.}\enspace\slant Let $E$ be reflexive and $A$ be a norm--closed linear monotone subspace of $E \times E^*$.   Then:
\smallbreak\noindent
{\rm(a)}\enspace $A$ is maximally monotone if, and only if, $A^*$ is monotone
\smallbreak
\noindent
{\rm(b)}\enspace $A$ is maximally monotone if, and only if, $A^*$ is maximally monotone.\endslant
\Proof These results follow from Theorem \PERPthm\ and the comments in Example \EEex.\qed
\defSection \APPsec
\medbreak
\centerline{\bf \APPsec.\quad The isometry lemma}
\medskip
\noindent
In this section, we make an observation about Banach SNL spaces that was not needed in the preceding sections.
\defLemma \FACTORlem
\medbreak
\noindent
{\bf Lemma \FACTORlem.}\enspace\slant Let $(B,L)$ be a Banach SNL space with Banach SNL dual $(B^*,\Lt)$.   Then $L$ is an isometry of $B$ into $B^*$.\endslant
\Proof Since $\|L\| \le 1$ and $\big\|\Lt\big\| \le 1$, for all $b \in B$,
$$\|b\| = \big\|\wh b\big\| = \big\|\Lt(Lb)\big\| \le \big\|\Lt\big\|\big\|Lb\big\| \le \big\|Lb\big\| \le \|L\|\|b\| \le \|b\|.$$
Thus\quad $\|Lb\| = \|b\|$.\qed
\medskip
Lemma \FACTORlem\ implies that if the surjectivity hypothesis, (\SURone), is satisfied then $B$ is what was called a \slant SSDB space\endslant\ in \cite\HBM, Definition 21.1\endcite, \cite\MLMAX, Section 1\endcite, \cite\SSDB, Definition 4.1\endcite\ and \cite\GMS, Section 3\endcite.   It also implies that $(\r,0)$ is a Banach SNL space without a Banach SNL dual.
\bigskip
\centerline{\bf References}
\medskip
\nmbr\ARENS
\item{[\ARENS]} R. Arens, \slant Operational calculus of linear relations\endslant, Pacific J. Math. {\bf 11} (1961) 9--23.
\nmbr\BB
\item{[\BB]} H. Brezis and F. E. Browder, {\sl  Linear maximal monotone operators and singular nonlinear integral equations of Hammerstein,  type}, Nonlinear analysis (collection of papers in honor of Erich H. Rothe),  pp. 31--42. Academic Press, New York, 1978.
\nmbr\BBWYLIN
\item{[\BBWYLIN]} H. Bauschke, J. M. Borwein, X. Wang and L. Yao \slant For maximally monotone linear relations, dense type, negative-infimum type, and Fitzpatrick-Phelps type all coincide with monotonicity of the adjoint\endslant, http://arxiv.org/abs/1103.6239v1, posted March 31, 2011.
\nmbr\BBWYFP
\item{[\BBWYFP]} H. Bauschke, J. M. Borwein, X. Wang and L. Yao \slant Every maximally monotone operator of Fitzpatrick-Phelps type is actually of dense type\endslant, http://arxiv.org/abs/1104.0750, posted April 5, 2011.
\nmbr\BBWYBB
\item{[\BBWYBB]} H. Bauschke, J. M. Borwein, X. Wang and L. Yao \slant The Brezis-Browder Theorem in a general Banach space\endslant, http://arxiv.org/abs/1110.5706v1, posted October 26, 2011.
\nmbr\BRON
\item{[\BRON]} A. Br\o ndsted and R.T. Rockafellar, \slant On the
Subdifferentiability of Convex Functions\endslant, Proc. Amer. Math.
Soc. {\bf 16}(1965), 605--611.
\nmbr\GMS
\item{[\GMS]}Y. Garc\'\i a Ramos, J. E. Mart\'\i nez-Legaz and S. Simons, \slant New results on $q$--positivity\endslant, http://arxiv.org/abs/1111.6094v1, posted November 25, 2011. 
\nmbr\ASD
\item{[\ASD]}M. Marques Alves and B. F. Svaiter, \slant On Gossez type (D) maximal monotone operators.\endslant, J. of Convex Anal., {\bf 17} (2010), 1077--1088.
\nmbr\MLMAX
\item{[\MLMAX]}J. E. Mart\'\i nez-Legaz, \slant On maximally $q$--positive sets\endslant, J. of Convex Anal., {\bf 16} (2009), 891--898.
\nmbr\FENCHEL
\item{[\FENCHEL]}R. T. Rockafellar, \slant Extension of Fenchel's duality theorem for convex functions\endslant, Duke Math. J. {\bf33} (1966), 81--89.
\nmbr\RANGE
\item{[\RANGE]}S. Simons, \slant The range of a monotone operator\endslant, J. Math. Anal. Appl. {\bf 199} (1996), 176--201.
\nmbr\PANDM
\item{[\PANDM]}-----, \slant Positive sets and Monotone sets\endslant, J. of Convex Anal., {\bf 14} (2007), 297--317.
\nmbr\HBM
\item{[\HBM]}-----, \slant From Hahn--Banach to monotonicity\endslant, 
Lecture Notes in Mathematics, {\bf 1693},\break second edition, (2008), Springer--Verlag.
\nmbr\SSDMON
\item{[\SSDMON]}-----, \slant Banach SSD spaces and classes of monotone sets\endslant, J. of Convex Anal., {\bf 18} (2011), 227--258.
\nmbr\SSDB
\item{[\SSDB]}-----, \slant SSDB spaces and maximal monotonicity\endslant, J. Glob. Optim., {\bf 50} (2011), 23--57.
\nmbr\YAO 
\item{[\YAO]} L. Yao, \slant The Brezis--Browder theorem revisited and properties of Fitzpatrick functions of order $n$\endslant, http://arxiv.org/abs/0905.4056v1, posted May 25, 2009.
\nmbr\ZBOOK
\item{[\ZBOOK]} C. Z\u{a}linescu, \slant Convex analysis in general vector spaces\endslant, (2002), World Scientific.
%
%
\Signoff
\bye